\documentclass[preprint,12pt]{elsarticle}
\usepackage{amssymb}
\usepackage{makeidx}
\usepackage{amsmath}
\usepackage{graphicx}
\usepackage{amsfonts}
\usepackage{pgf}
\usepackage{tikz}
\usetikzlibrary{external}
\tikzset{external/up to date check=simple}
\usepackage{xifthen}
\newboolean{pedro}
\setboolean{pedro}{true}
\newcommand{\pedro}{\ifthenelse{\boolean{pedro}}{\color{black}
    \setboolean{pedro}{false}}{\color{black}\setboolean{pedro}{true}}}
\newcounter{margin}

\newboolean{pedrom}
\setboolean{pedrom}{true}
\newcommand{\pedrom}{\ifthenelse{\boolean{pedrom}}{\color{black}
    \setboolean{pedrom}{false}}{\color{black}\setboolean{pedrom}{true}}}

\title{Optimal Control of Counter-Terrorism Tactics}
\author{L. Bay\'on\corref{cor1}}
\author{P. Fortuny Ayuso\corref{cor2}}
\author{P. J. Garc\'ia-Nieto\corref{cor2}}
\author{J. M. Grau\corref{cor2}}
\author{M. M. Ruiz\corref{cor2}}

\address{Department of Mathematics, University of Oviedo, E.P.I. Campus of Viesques, Gij\'{o}n, 33203, Spain.}
\cortext[cor1]{Corresponding author: bayon@uniovi.es}

\newproof{pf}{Proof}

\journal{Applied Mathematics and Computation}

\begin{document}

\begin{abstract}
  This paper presents an optimal control problem to analyze the efficacy of
  counter-terrorism tactics. We present an algorithm that efficiently combines
  the Minimum Principle of Pontryagin, the shooting method and the cyclic
  descent of coordinates. We also present a result that allows us to know a
  priori the steady state \pedro{}solutions\pedro{}. Using this technique we are
  able to choose parameters that reach a specific solution, \pedro{}of which
  there are two\pedro{}. Numerical examples are presented to illustrate the
  possibilities of the method. Finally, we study the sufficient conditions for
  optimality and suggest an improvement on the functional which also guarantees
  local optimality.
\end{abstract}


\maketitle

\section{Introduction}

Modeling a ``stock'' of terrorists, is not common, but has precedents,
especially after September 11, 2001 \cite{Heymann 2003}. In this sense
\cite{Keohane 2003} presents an intelligent ecological metaphor to analyze
actions by Governments and citizens against terror. In \cite{Castillo 2003} a
model for the transmission dynamics of extreme ideologies in vulnerable
populations is presented. In \cite{Kaplan 2005} the authors propose a
terror-stock model that treats the suicide bombing attacks in Israel. In other
countries like, for example, Spain or Ireland, the problem has also been analyzed.

Several papers develop dynamical models of terrorism. In \cite{Udwadia 2006}
the authors incorporate the effects of both military/police and
nonviolent/persuasive intervention to reduce the terrorist population. This
idea is widely developed in \cite{Caulkins 2008} where the controls are two
types of counter-terror tactics: ``water'' and ``fire'', which is the model we
shall consider in this paper. The effect of intelligence (water tactics) in
counter-terrorism is analyzed also in\ \cite{Kress 2009}. Nowadays, it is
agreed that counter-terrorism policies have the potential to generate positive
support for terrorism \cite{Faria 2012}. Recently, in \cite{Seidel 2016} a
model with two-states (undetected and detected terrorists) and only one
control variable (the number of undercover intelligence agents) is considered.

In this context we present in this work a new approach to analyze the efficacy
of counter-terrorism tactics. We state an optimal control problem that
attempts to minimize the total cost of terrorism. An excellent summary of
optimal control application in these issues can be consulted in
\cite{Grass 2008} and its economic implications in \cite{Schneider 2015}.

The optimization criterion is to minimize the discounted damages created by
terror attacks plus the costs of counter-terror efforts. The underlying
mathematical problem is complicated. It constitutes a multi\pedrom{}-control\pedrom{},
constrained problem where the optimization interval is infinite. An important
feature is that the time $t$ is not explicitly present in the problem (hence, it
is a time-autonomous problem), except in the discount factor. Using Pontryagin's
Minimum Principle, the shooting method and the cyclic descent of coordinates we
develop an optimization algorithm. We also present a method (based upon
\cite{Tsur 2001}) for computing the optimal \pedro{}steady-states\pedro{} in
multi\pedrom{}-control\pedrom{}, infinite-horizon, autonomous models. This method does not require
the solution of the dynamic optimization problem. Using it, we can choose
parameters that reach a desirable steady-state solution. \pedro{}The problem
presents two steady-states, albeit one of them in a region where it becomes
effectively one-dimensional. We focus mainly on the multi\pedrom{}-control\pedrom{}
problem.\pedro{}

The paper is organized as follows. Section 2 presents the mathematical model.
The optimization algorithm is developed in Section 3, and the method for
computing the optimal \pedro{}steady-states\pedro{} is analyzed in Section
4. Section 5 presents several numerical examples which illustrate the
performance of the algorithm under different conditions.  In Section 6 we
discuss Arrow's sufficient conditions for optimality in our problem. We also
suggest an improvement on the functional in which the cost function is convex in
the number of terrorists (due to the value added by information sharing,
interactions, etc.) and show how the solution found by our method in this case
satisfies the sufficient conditions locally.  Finally, the main conclusions of
our work are discussed in Section 7.

\section{Mathematical Model}

We use the excellent model provided by \cite{Caulkins 2008},
which classifies counter-terrorism tactics into two categories:

\begin{itemize}
\item ``Fire'' strategies are tactics that involve significant collateral
damage. They  include, for example, the killing of terrorists through drones,
the use of  indiscriminate checkpoints or the aggressive blockade of roads.

\item ``Water'' strategies, on the other hand, are counter-measures that do
not affect  innocent people, like intelligence arrests against suspect individuals.
\end{itemize}

The fire and water strategies will be denoted by the control variables $v(t)$
and $u(t)$, respectively. Both controls have their advantages, and their
drawbacks. For example $v(t)$ has the direct benefit of eliminating current
terrorists but the undesirable indirect effect of stimulating recruitment
rates (and the possible harm to innocent bystanders). On the other hand,
$u(t)$ is more expensive and more difficult to be applied than $v(t)$.

The strength or size of the terrorists is represented by the state variable
$x(t)$. This includes not only the number of active terrorists, but also the
organization's total resources including financial ones, weapons, etc.
\cite{Keohane 2003}. Its value changes over time and we distinguish two
inflows and three outflows in it:%
\begin{equation}
\dot{x}=\tau+I(v,x)-O_{1}(x)-O_{2}(u,x)-O_{3}(v,x) \label{Dynamic}%
\end{equation}
We include first of all a term $\tau$, accounting
for a small constant recruitment rate. Second, following \cite{Castillo 2003},
the model considers that new terrorists are recruited by existing terrorists.
So the inflow $I(v,x)$ is increasing in proportion to the current number of
terrorists $x$. But this growth is bounded and should also slow down.
Moreover, the aggressive control $v$, also increases recruitment. In summary
the form of the model is:%
\begin{equation}
I(v,x)=(1+\rho v)kx^{\alpha} \label{Inflow}%
\end{equation}
with $\tau,\rho\geq0,$ $k>0$ and $0\leq\alpha\leq1$.

On the other hand, we consider three outflows: The first one, $O_{1}(x)$,
represents the rate at which people leave the organization by several reasons
not related with the controls. This natural outflow is assumed linear in $x$:
\begin{equation}
O_{1}(x)=\mu x \label{O1}%
\end{equation}
with $\mu>0.$ The second outflow, $O_{2}(u,x)$, reflects the effect of water
strategies. This outflow is assumed to be concave in $x$ because there is a
limited number of units that conduct water operations:%
\begin{equation}
O_{2}(u,x)=\beta(u)x^{\theta} \label{O2}%
\end{equation}
with $\theta\leq1.$ The third outflow $O_{3}(v,x)$ is due to fire strategies.
This is modeled as linear in $x$, because the methods are perceived to be
``direct attack'':\
\begin{equation}
O_{3}(v,x)=\gamma(v)x \label{O3}%
\end{equation}
The functions $\beta(u)$ and $\gamma(v)$ should be concave; Caulkins
\cite{Caulkins 2008} uses the same functional form for both: a logarithmic
function. The water function is pre-multiplied by a constant $\beta$ smaller
than the corresponding constant $\gamma$ for fire operations. These two
constants reflect the ``efficiency'' of the two types of operations.

Finally, the costs of terrorism are assumed to be linear in the number of
terrorists, that is, of the form $cx$. We also model the cost control function
as separable, and the costs of employing the water and fire strategies are
modeled as quadratic. Over a infinite planning horizon, the objective is to
minimize the sum of both costs (terrorism and counter-terror operations). We
also assume that outcomes are discounted by a constant rate $r$. In brief, the
control problem we pose can be written as:%
\begin{align}
\underset{u,v\geq0}{\min}J  &  =\underset{u,v\geq0}{\min}\int_{0}^{\infty
}(cx+u^{2}+v^{2})e^{-rt}dt\label{OCP}\\
\dot{x}  &  =\tau+(1+\rho v)kx^{\alpha}-\mu x-\beta\ln(1+u)x^{\theta}%
-\gamma\ln(1+v)x;\text{ \ }x(0)=x_{0}\nonumber\\
u(t)  &  \geq0;\text{ \ }v(t)\geq0\nonumber
\end{align}
where $x_{0}$ is the initial stock level and we impose also control constraints.

\section{Optimization Algorithm}

The above problem (\ref{OCP}), is an Optimal Control Problem (OCP) where the
total costs have to be minimized, given the state dynamics and the control
constraints. Denoting $\mathbf{u}(t)=(u(t),v(t))=(u_{1}(t),u_{2}(t))$ we have:%
\begin{equation}
\underset{\mathbf{u}(t)}{\min}J=\int_{0}^{\infty}F\left(  t,x(t),\mathbf{u}%
(t)\right)  dt\label{Funct}%
\end{equation}
subject to satisfying:
\begin{align}
\dot{x}(t) &  =f\left(  t,x(t),\mathbf{u}(t)\right)  ,\text{ }0\leq
t<\infty;\text{ \ }x(0)=x_{0}\label{state eq}\\
\mathbf{u}(t) &  \in\mathbf{U}(t),\text{ }0\leq t<\infty\label{bounds}%
\end{align}
The problem presents several noteworthy features. First, the optimization
interval is infinite. Second, the time $t$ is not explicitly present in the
problem (time-autonomous problem), except in the discount factor. Third, we
impose constraints on the control and, fourth, it constitutes a
multi\pedrom{}-control\pedrom{} problem.

\subsection{Multi-control problem}

To solve the multi\pedrom{}-control\pedrom{} variational problem, we propose a numerical
algorithm which uses a particular strategy related to the cyclic coordinate
descent (CCD) method \cite{Luo 1992}. The classic CCD method minimizes a
function of $n$ variables cyclically with respect to the coordinates. With our
method, the problem can be solved like a sequence of problems whose error
functional converges to zero. The algorithm (with $i=1,2$) carries out several
iterations and at each $j$-\textit{th }iteration it calculates $2$ stages, one
for each $i$. At each stage, it computes the optimal of $u_{i}(t)$, assuming
the other variable is fixed.

Beginning with some admissible $\mathbf{u}^{0}$, we construct a sequence of
$(\mathbf{u}^{j}\mathbf{)}$ and the algorithm will search:
\begin{equation}
\underset{j\rightarrow\infty}{\lim}\mathbf{u}^{j} \label{conv}%
\end{equation}
It is easy to justify the convergence of the algorithm taking into account
Zangwill's global convergence Theorem \cite{Zangwil 1969}.

\pedro{}\paragraph{The whole problem is not multi-control} We
remark here, in passing, that, as we shall see later, there is a value $x^S$,
which we call a \emph{switching point}, such that the optimal solution
$(x^{\ast}(t),u^{\ast}(t),v^{\ast}(t))$ satisfies $v^{\ast}(t)=0$ for all $t$
such that $x^{\ast}(t)\leq x^{S}$. This switching point is independent of the
solution; therefore, for $x\leq x^{S}$, the problem becomes essentially
uni\pedrom{}-control\pedrom{}. However, for $x>x^{S}$ the multi\pedrom{}-control\pedrom{} techniques are
required.\pedro{}

\pedro{}\subsection{The one-dimensional reduced problem}\pedro{}

Based upon the previous statement, we present now the solution for the
one-dimensional case, using Pontryagin's Minimum Principle (PMP) (see, for
example, \cite{Pontryagin 1987}, \cite{Clarke 1983}, \cite{Chiang 2000}).
\begin{align}
\underset{u(t)}{\min}J  & =\int_{0}^{\infty}F\left(  t,x(t),u(t)\right)
dt\label{OCP 1}\\
\dot{x}(t)  & =f\left(  t,x(t),u(t)\right)  ,\text{ }0\leq t<\infty;\text{
\ }x(0)=x_{0}\nonumber\\
u(t)  & \geq0,\text{ }0\leq t<\infty\nonumber
\end{align}
Our integrand takes the form:%
\begin{equation}
F(t,x(t),u(t))=G(t,x(t),u(t))e^{-rt}%
\end{equation}
where $r$ is the positive rate of discount, and $G$ is a function bounded from
above. Under these conditions, the integral is found to be convergent for each
admissible control. Let $H$ be the associated Hamiltonian:
\begin{equation}
H(t,x,u,\lambda)=\pedrom{}\lambda^0\pedrom{}F\left(  t,x,u\right)  +\lambda\cdot f\left(  t,x,u\right)
\label{Ham}%
\end{equation}
where $\lambda$ is the co-state variable. Using PMP, the optimal solution can
be obtained from a two-point boundary value problem. In order for $u^{\ast}\in
U$ to be optimal, there must exist a function $\lambda(t)$ such that, for
almost every $t\in\lbrack0,\infty)$:
\begin{align}
\dot{x} &  =H_{\mathbf{\lambda}}=f;\text{ \ }x(0)=x_{0}\\
\dot{\lambda} &  =-H_{\mathbf{x}};\text{ \ }\lim_{t\rightarrow\infty}
\pedrom{}H(t)\pedrom{}=0\label{lambda p}\\
H(t,x,u^{\ast},\lambda) &  =\underset{u(t)\in U}{\min}H(t,x,u,\lambda
)\label{Min H}%
\end{align}
\pedrom{}An elementary argument shows that, in the problem under consideration, any optimal path $(x^{\ast}(t), u^{\ast}(t), v^{\ast}(t))$ satisfies that $0<u^{\ast}(t)<M$ for some $M>0$. This, together with the transversality condition  $\lim_{t\rightarrow \infty}H(t)=0$ implies that both $\lambda^{0}\neq 0$ and $\lim_{t\rightarrow\infty}\lambda(t)=0$ (which is sometimes stated as a transversality condition). The latter is the property we shall use later on in order to compute an optimal path. Also, as $\lambda^0\neq 0$, we normalize it to $1$.\pedrom{}

\pedrom{}Notice that, frequently, the condition
\begin{equation}
\lim_{t\rightarrow\infty}\lambda(t)=0\label{TC}%
\end{equation}
is stated as a necessary condition for optimality in infinite-time problems. In  \cite{Benveniste 1982}, it is proved that, under certain hypotheses, for problems with time discount, it is necessary. However, those hypotheses can be (as in this case) not easily checked. Notice also, that there are examples of infinite-time problems with discount in which $\lambda(t)$ does not converge to $0$ \cite{aseev}.
\pedrom{}

Due to the non-linearity of the system dynamics, the optimal solution can only
be computed numerically. In this paper we propose an efficient method which
adapts the shooting method, Euler's method, and numerical integration. All the
calculations are carried out in the Mathematica environment.

We will denote by $\mathbb{Y}_{x}(t)$ the function:
\begin{equation}
\mathbb{Y}_{x}(t)=-\frac{F_{u}}{f_{u}}\cdot e^{\int_{0}^{t}f_{x}ds}+\int
_{0}^{t}F_{x}\cdot e^{\int_{0}^{s}f_{x}dz}ds\label{COORD F}%
\end{equation}
The algorithm is based upon the following theorem.

\textbf{Theorem 1. A necessary maximum condition}

\textit{Let }$u^{\ast}$ \textit{be the optimal control, let} $x^{\ast}%
\in\widehat{C}^{1}$ \textit{be a solution of the above problem (\ref{OCP 1}).
Then there exists a constant }$K\in\mathbb{R}^{+}$\textit{\ such that:}%
\begin{equation}%
\begin{tabular}
[c]{lll}%
\textit{If }$u^{\ast}(t)>0$ & $\Longrightarrow$ & $\mathbb{Y}_{x^{\ast}}(t)=K
$\\
\textit{If }$u^{\ast}(t)=0$ & $\Longrightarrow$ & $\mathbb{Y}_{x^{\ast}%
}(t)\geq K$%
\end{tabular}
\end{equation}

\textbf{Proof.}

In virtue of PMP there exists a piece-wise $C^{1}$\ function $\lambda$\ that
satisfies the linear differential equation:
\begin{equation}
\dot{\lambda}(t)=-H_{x}=-F_{x}-\lambda(t)\cdot f_{x}%
\end{equation}
Denoting $K=\lambda(0)$, we have:
\begin{equation}
\lambda(t)=\left[  K-\int_{0}^{t}F_{x}e^{\int_{0}^{s}f_{x}dz}ds\right]
e^{-\int_{0}^{t}f_{x}ds}\label{lambda}%
\end{equation}
For each $t$, $u(t)$ minimizes $H$. Using the Kuhn-Tucker Theorem, there
exists for each $t$ a real non negative number, $\delta(t)$, such that $u(t)$
is a critical point of the augmented Hamiltonian:
\begin{equation}
\mathbb{H}(u(t))=F(t)+\lambda(t)\cdot f(t)+\delta(t)\cdot(-u(t))\label{K-T}%
\end{equation}
If $u(t)>0$, then $\delta(t)=0$. In this case:
\begin{equation}
F_{u}+\lambda(t)\cdot f_{u}=0\label{der}%
\end{equation}
\pedrom{}Notice that $f_u\neq 0$ in our problem, as $f_u=-\beta(1+u)^{-1}x^{\theta}$. Using this, \pedrom{}from (\ref{lambda}) and (\ref{der}), we obtain:
\begin{equation}
K=-\frac{F_{u}}{f_{u}}\cdot e^{\int_{0}^{t}f_{x}ds}+\int_{0}^{t}F_{x}\cdot
e^{\int_{0}^{s}f_{x}dz}ds\label{EQ COOR}%
\end{equation}
and the following relation holds:
\begin{equation}
\mathbb{Y}_{x}(t)=K\label{T1}%
\end{equation}
If $u(t)=0$, then $\delta(t)\geq0$. By a similar argument and bearing in mind
that $f_{u}\leq0$, we get that
\begin{equation}
\mathbb{Y}_{x}(t)\geq K\label{T3}%
\end{equation}

$\hfill\square$

We have performed the numerical (and approximate) construction of $x^{\ast}$
using a discretized version of equation (\ref{EQ COOR}). For each $K$, we
construct $x^{\ast}$ using (\ref{T1}) and (\ref{T3}) to impose the constraint.

The calculation of the optimal $K$ has been achieved by means of an adaptation
of the shooting method. Varying the constant, $K$, we search for the extremal
that fulfills the transversality condition (\ref{TC}). Using the secant method
(starting out from two values $K_{\min}$ and $K_{\max}$), our algorithm
converges satisfactorily (see Section 5).

\pedrom{}Finally: an adapted version of the CCD method for functionals (whose
detailed description can be found in \cite{bolza}, applied to hydrotermal
problems), has been used in order to solve the multi-control problem. The
classic CCD method minimizes a function of n variables cyclically with respect
to the coordinate variables. With our method, the problem is solved as a
sequence of problems whose error functional converges to zero.\pedrom{}

%Finally, we have used an adapted version of the CCD method for functionals, to
%solve the multi\pedrom{}-control\pedrom{} problem.

\pedro{}\section{Steady-state solutions}\pedro{}

In \cite{Tsur 2001} a method for computing the optimal steady-state in
infinite-horizon one-dimensional problems is presented which does not require
the solution of the dynamic optimization problem, in which the bounds $U(t)$
do not play any role. Tsur considers a one-dimensional version of our problem:%
\begin{align}
\underset{\mathbf{u}(t)}{\min}J &  =\int_{0}^{\infty}G(x(t),\mathbf{u}%
(t))e^{-rt}dt\\
\dot{x}(t) &  =f\left(  x(t),\mathbf{u}(t)\right)  ,\text{ \ }x(0)=x_{0}%
\end{align}
We propose another adaptation of the CCD method. Beginning with some
admissible $\mathbf{u}^{0}$, we construct a sequence $(\mathbf{u}^{j})$ and at
each stage, we compute the optimal steady state of $u_{i}(t)$, assuming the
other variable is fixed. At each $i$-th stage (for $i=1,2$), we consider that
the steady-state solution is $u_{i}=R_{i}(x)$ and we define the
\textit{evolution function}:%
\begin{equation}
L_{i}(x)=r\left(  \frac{G_{u_{i}}(x,R_{i}(x))}{f_{u_{i}}(x,R_{i}(x))}+\dot
{W}_{i}(x)\right)  \label{L}%
\end{equation}
with:%
\begin{equation}
W_{i}(x)=\frac{1}{r}G(x,R_{i}(x))\label{W}%
\end{equation}
A necessary condition for \pedro{}an\pedro{} optimal steady state $x_{s}$ is:
\begin{equation}
L_{i}(x_{s})=0\label{CN}%
\end{equation}
The algorithm shows a fast convergence to the optimal values of $u_{i}(t)$ and
also gives the unique value of $x_{s}\pedro{}>x^S\pedro{}$. It is noteworthy
that for the dynamic equation model presented in (\ref{OCP}) with
$\mathbf{u}(t)=(u(t),v(t))=(u_{1}%
(t),u_{2}(t))$:%
\begin{equation}
\dot{x}=\tau+(1+\rho v)kx^{\alpha}-\mu x-\beta\ln(1+u)x^{\theta}-\gamma
\ln(1+v)x
\end{equation}
the function $u=R_{1}(x),$ obtained by imposing $\dot{x}=0$, is:%
\begin{equation}
R_{1}(x)=-1+\exp\left(  \frac{x^{-\theta}(\tau+(1+\rho v)kx^{\alpha}-\mu
x-\gamma\ln(1+v)x)}{\beta}\right)
\end{equation}
Nevertheless, the value of $v=R_{2}(x),$ is not so easy to obtain. It turns
out to be
\begin{equation}
R_{2}(x)=-1-\frac{x^{1-\alpha}\gamma}{k\rho}W\left(  -\frac{1}{\gamma}%
\exp\left(  \dfrac{\tau-\mu x-kx^{\alpha}(\rho-1)}{x\gamma}\right)
k(1+u)^{\dfrac{-x^{\theta-1}\beta}{\gamma}}x^{\alpha-1}\rho\right)
\end{equation}
where $W$ is the Lambert $W$-function. Remember that $W(z)$ is a set of
functions which are the branches of the inverse of the function:
\begin{equation}
z=f(W)=We^{W}%
\end{equation}
where $W$ is a complex variable. In this work we are interested in real-valued
$W(x)$, which, adding the condition $W(x)\geq-1$, gives a single-valued
function $W_{0}(x)$: the principal branch of the $W$-function. We refer the
reader to \cite{Corless 1996} for a survey on existing results on this function.

In the next section we will see the behavior of this algorithm and its
application to calculate the optimal steady-state of our problem, \pedro{}for
$x>x^S$\pedro{}.

\pedro{}\paragraph{Steady-states before the switching point}
As remarked above, the problem becomes effectively
one-dimensional when $x\leq x^S$: namely $v=0$ for $x\leq x^S$. This simplifies
the problem in this case, obviously. \pedrom{}As we shall see\pedrom{}, Tsur's method provides the
existence of two steady states in $x\in[0,x^S]$: one of them stable and the
other one unstable, without recoursing to multi\pedrom{}-control\pedrom{}
techniques. Therefore, the optimization problem presents two stable steady
states, one to the left of the switching point and another one to the
right.\pedro{}

\section{Numerical Examples}

\subsection{Base Case}

We examine now the behavior of our approach in several examples. For the sake
of comparison, we use the (carefully chosen) parameters used in \cite{Caulkins
2008}. The discount rate is a typical $r=0.05$. The outflow rate is assumed to
be $5\%$ \pedrom{}and the constant\pedrom{} inflow rate term is small $\tau=10^{-5}$. The
parameter $k$ is chosen such that the steady state is normalized to $1$ in the
absence of counter-terrorism tactics, and neglecting $\tau$. This way, $x$ is
measured as a percentage of the steady-state size of the terrorist
organization when the government does not use counter-terror operations. The
uncontrolled dynamics is given by:%
\begin{equation}
\dot{x}=kx^{\alpha}-\mu x
\end{equation}
Integrating this Bernoulli differential equation, and computing the limit when
$t\rightarrow\infty$ we have:%
\begin{equation}
x=\sqrt[1-\alpha]{\dfrac{k}{\mu}}%
\end{equation}
We see that the normalization of $x$ leads to $k=\mu.$ The influence of $x$ on
recruitment, $\alpha$, is a value between $0$ and $1$. The efficiency
parameters $\gamma$ and $\beta$ are chosen assuming fire strategies are
approximately $10$ times more powerful than water strategies for the maximum
size of $x$. Nevertheless, for small values of $x$, fire and water strategies
can be equally effective. In this case, we choose the parameter $\theta$ such
that (for small $x\simeq1/10$):
\begin{equation}
\beta\ln(1+u)x^{\theta}=\gamma\ln(1+v)x\rightarrow\frac{\gamma}{\beta
}=x^{1-\theta}=1
\end{equation}
The parameter $\rho$ measures the influence of fire strategies on recruitment.
It is chosen so that the product $\rho v(t)$ represents the $20-30\%$ of the
recruitment rate. Finally, the costs $c$ per terrorist are assumed to be $1$
in the base case. Table 1 summarizes the values of the parameters for the
model.%
\begin{gather*}
\text{Table 1. Parameters for the fire-and-water model \cite{Caulkins 2008}%
.}\\%
\begin{tabular}
[c]{lll}%
\textbf{Par.} & \textbf{Description} & \textbf{Value}\\\hline
$r$ & Discount rate & 0.05\\
$c$ & Costs per unit of terrorists & 1\\
$\tau$ & Constant inflow & 10$^{-5}$\\
$\rho$ & Contribution of fire control to recruitment & 1\\
$k$ & Normalization factor & 0.05\\
$\alpha$ & Influence of actual state on recruitment & 0.75\\
$\mu$ & \textquotedblleft Natural\textquotedblright\ per capita outflow &
0.05\\
$\beta$ & Efficiency of water operations & 0.01\\
$\theta$ & Diminishing returns from water operations & 0.1\\
$\gamma$ & Maximum efficiency of fire operations & 0.1
\end{tabular}
\end{gather*}
\pedro{}\paragraph{Switching point} From the equations found in
\cite{Caulkins 2008} (which we omit for brevity), the value of the optimal
control $v^{\ast}(t)$ is, without active restrictions (and using our notation,
in which we do not remove the discount factor)
\begin{equation}\label{eq:vopt}
  v^{\ast}(t) = \frac{1}{4}\left(
    -\rho k x(t)^{\alpha}\lambda(t)e^{-rt}-2+
    \sqrt{(\rho k x(t)^{\alpha}\lambda(t)e^{-rt}-2)^2 +
      8 \lambda(t)e^{-rt}\nu x(t) }
  \right).
\end{equation}
Let us remark that the equation given for $v^{\ast}$ in Equation (3) in
\cite{Caulkins 2008} lacks a factor $\lambda$ in the first term out of the square root, which prevents the authors from exactly computing the switching
point. Equation \eqref{eq:vopt} has the form (omitting the factor $1/4$)
\begin{equation}
  v^{\ast}(t)=-(a+b)+\sqrt{(a-b)^2+8c} = -(a+b)+\sqrt{(a+b)^2-4ab+8c}.
\end{equation}
Rewriting \eqref{eq:vopt} in this form, we get
\begin{equation}
  v^{\ast}(t) = -(2+\rho k x(t)^{\alpha}\lambda(t)e^{-rt}) +
  \sqrt{(2+\rho k x(t)^{\alpha}\lambda(t)e^{-rt})^2-8 \lambda(t)e^{-rt}(\rho k x(t)^{\alpha} -  \nu x(t))}
\end{equation}
so that (as $\lambda(t)> 0$) whenever $\rho k x(t)^{\alpha}>\nu x(t)$, the
optimal solution has $v^{\ast}(t)< 0$, which is not possible (as $v(t)$
denotes a positive \emph{action}). Hence, if
\begin{equation}
  \rho k x(t)^{\alpha} > \nu x(t)
\end{equation}
the restriction $v(t)=0$ is activated. This is the same as
\begin{equation}
  x(t)^{\alpha}(\rho k - \nu x(t)^{1-\alpha})>0
\end{equation}
with $x(t)\geq 0$. That is, whenever
\begin{equation}
  x(t) < x^S=\left(\frac{\rho k}{\nu}\right)^{\frac{1}{1-\alpha}}
\end{equation}
one has to impose $v^{\ast}(t)=0$. For the given values of the parameters, we
obtain that this switching point is $x^S=0.0625$ \emph{exactly}.  \pedro{}

\pedro{}\subsection{High-value steady state}\pedro{}
We implemented a Mathematica program to apply the results obtained in this
paper. \pedro{}As noted in the previous section, the steady-state
solutions can be computed \emph{a priori}. We focus on the high value of the
steady state, $x_s$, as this is the one in which the system is truly
multi\pedrom{}-control\pedrom{}. As there are no more steady states with $x>x_s$, we shall
compute it starting from $x(0)>x_s$.\pedro{} The algorithm, adapting the CCD
method, solves equation (\ref{CN}) in $x_{s}$. Starting from two initial values
$u^{0}(t)=v^{0}(t)=0$, we obtain the following steady-state values:
\begin{align}
x_{s} &  \simeq 0.61773\\
u_{s} &  \simeq 0.06834\\
v_{s} &  \simeq 0.14605
\end{align}
They were also obtained in \cite{Caulkins 2008} (but only approximately: he
gives $x_{s}\simeq0.62$). However, it should be recalled that the
determination of this value is as important as the computation of the dynamics
of the process leading towards that steady-state.

In order to compute the dynamic solution, we use the optimization algorithm
\pedrom{}based on Theorem 1\pedrom{}. The following results were obtained using a discretization
in $t$ with $250$ sub-intervals. Starting with $u^{0}(t)=v^{0}(t)=0$ and two
initial values of $K$ for each control, namely $K_{\min}$ and $K_{\max}$, we
compute the discretized solution $x^{\ast}$ step by step with:
\begin{equation}
\mathbb{Y}_{x}(t)=K\text{ or }\mathbb{Y}_{x}(t)\geq K
\end{equation}
Using the secant method, we achieve the prescribed tolerance ($10^{-5}$) in
the transversality condition
\begin{equation}
\lim_{t\rightarrow\infty}\lambda(t)=0
\end{equation}
in only $9$ iterations for $u(t)$ (see Fig. 1(a)), and in $50$ iterations for
$v(t)$ (see Fig. 1(b)). Once these values are computed, the algorithm
alternates the two controls in a cyclic way, retro-feeding the values
$u^{j}(t),v^{j}(t)$ computed in each iteration. The CCD algorithm shows a fast
convergence: in just $4$ iterations the tolerance imposed for the error is
reached: in our case, the difference between two consecutive values of $K$ to
be less than $10^{-3}$ (see Fig. 2).
\begin{gather*}
\includegraphics[height=1.075in]{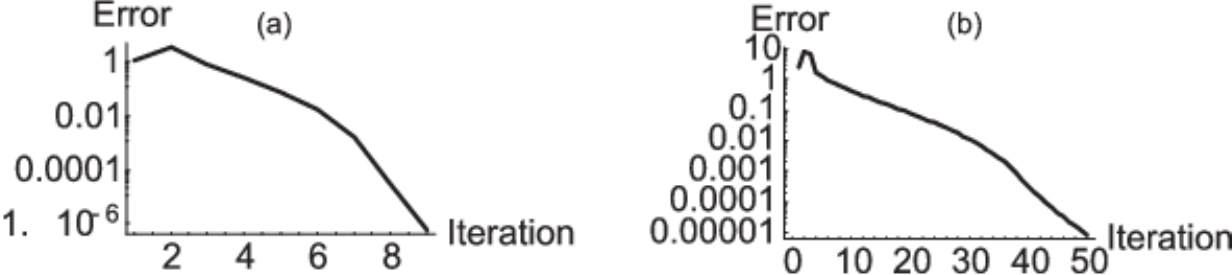}\\
  \text{Figure 1. Convergence of the algorithm for each control. }%
\end{gather*}
Our results show the double convergence of the algorithm: on one hand, of each
control to the transversality condition (\ref{TC}) and, on the other, of the
CCD method. We note this latter result especially, as it not only shows the
fast convergence but it also opens up the possibility of applying our method
to problems with more dimensions.
\begin{gather*}%
%TCIMACRO{\FRAME{itbpF}{2.2866in}{1.075in}{0in}{}{}{fig-2.pdf}%
%width 2.2866in;  height 1.075in;  depth 0in;  original-width 6.0727in;
%original-height 2.8392in;  cropleft "0";  croptop "1";  cropright "1";
%filename 'Fig-2.pdf';file-properties "NPEU";}%
%}}%
%BeginExpansion
{\includegraphics[
height=1.075in,
width=2.2866in
]%
{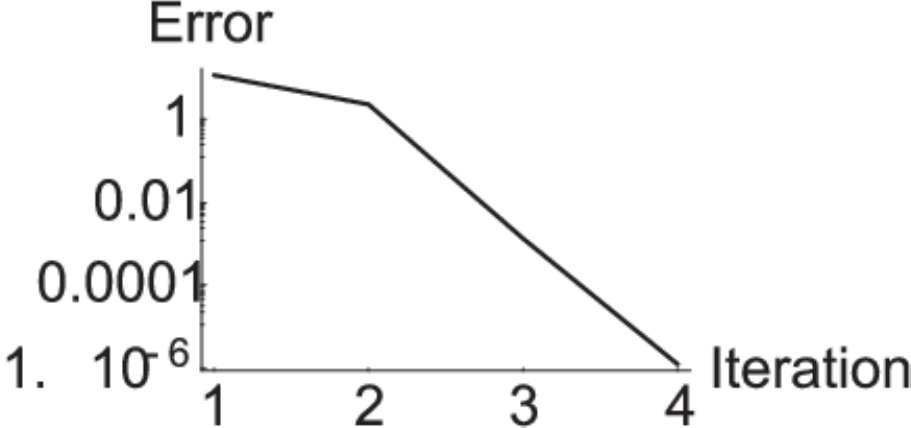}%
}%
%EndExpansion
\\
\text{Figure 2. Convergence of the CCD algorithm. }%
\end{gather*}
The optimal dynamical evolution of the strength or size of the terrorists,
$x(t)$, is shown in Fig. 3 and of the fire\pedrom{}-control\pedrom{} $v(t)$, is shown in Fig. 4.
The optimal path has been obtained with an initial value $x(0)=0.95$. The plot
of the water\pedrom{}-control\pedrom{} $u(t)$ is not shown because its value is practically
constant in time. As a matter of fact, it starts at $0.066$ for $t=0$ and
increases slowly towards its steady state value: $0.068$.
\[%
\begin{tabular}
[c]{ccc}%
$%
%TCIMACRO{\FRAME{itbpF}{2.0418in}{1.1519in}{0in}{}{}{fig-3.pdf}%
%width 2.0418in;  height 1.1519in;  depth 0in;  original-width 5.9352in;
%original-height 3.3131in;  cropleft "0";  croptop "1";  cropright "1";
%filename 'Fig-3.pdf';file-properties "NPEU";}%
%}}%
%BeginExpansion
{\includegraphics[
width=2in
]%
{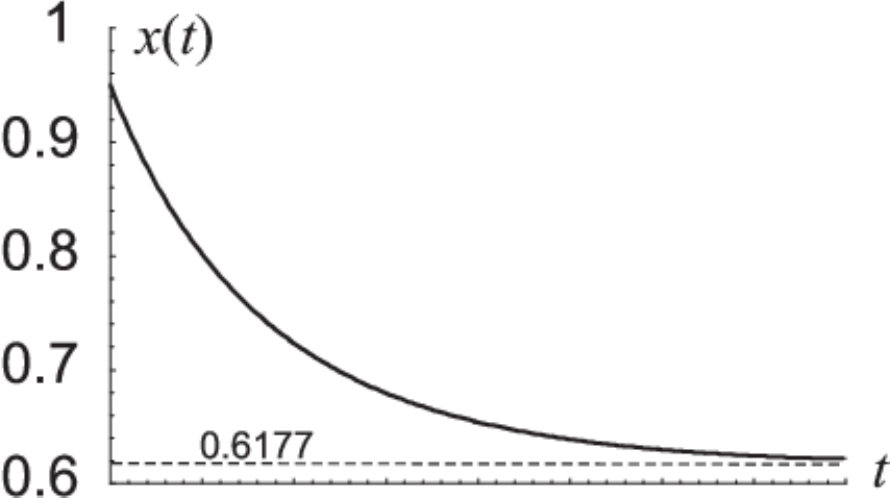}%
}%
%EndExpansion
$ &
%TCIMACRO{\FRAME{itbpF}{2.1499in}{1.1519in}{0in}{}{}{fig-4.pdf}%
%width 2.1499in;  height 1.1519in;  depth 0in;  original-width 7.3059in;
%original-height 3.8925in;  cropleft "0";  croptop "1";  cropright "1";
%filename 'Fig-4.pdf';file-properties "NPEU";}%
%}}%
%BeginExpansion
{\includegraphics[
width=2in
]%
{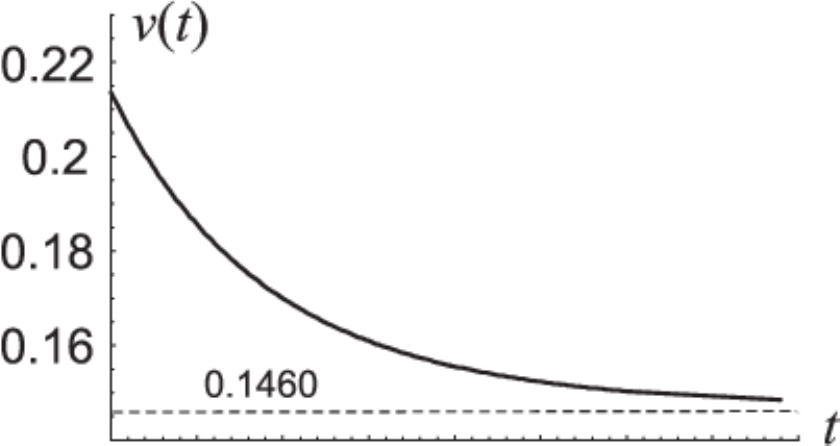}%
}%
%EndExpansion
& $%
%TCIMACRO{\FRAME{itbpF}{2.0193in}{1.1519in}{0in}{}{}{fig-5.pdf}%
%width 2.0193in;  height 1.1519in;  depth 0in;  original-width 6.0424in;
%original-height 3.4117in;  cropleft "0";  croptop "1";  cropright "1";
%filename 'Fig-5.pdf';file-properties "NPEU";}%
%}}%
%BeginExpansion
{\includegraphics[
width=2in
]%
{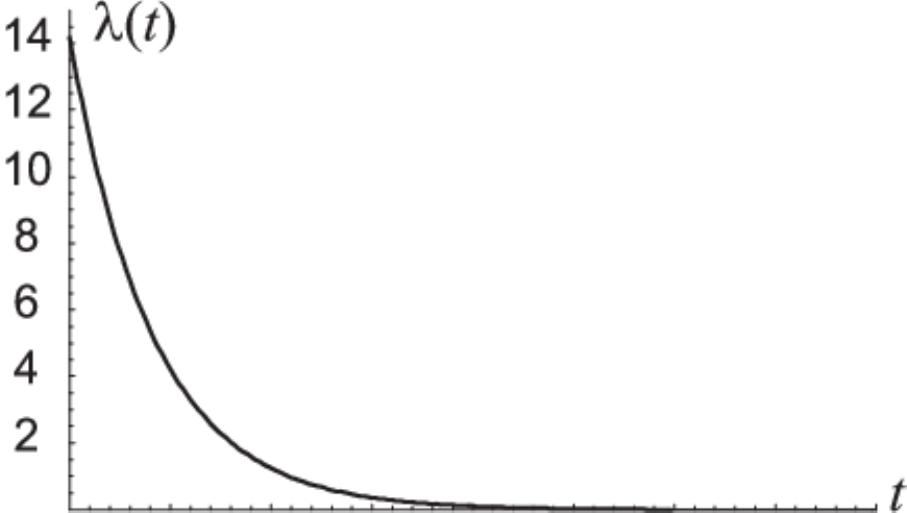}%
}%
%EndExpansion
$\\
$\text{Figure 3. Optimal }x(t)\text{.}$ & $\text{Figure 4. Optimal
 }v(t)\text{.}$ & $\text{Figure 5. Optimal }\lambda(t)\text{.
}$%
\end{tabular}
\]
The CPU time required by the program was $121$ $\sec$ on a personal computer
(Intel Core 2/$2.66$ GHz). The optimal value of $K$ is $14.16651221278533$.
The evolution of $\lambda(t)$ and its asymptotic behavior towards zero, can
be seen in Fig. 5.

\subsection{Low-value steady state}
\pedro{}
Using
Tsur's method straightforwardly, we can compute the steady states for $x<x^S$,
as the problem is one-dimensional ($v(t)=0$ for these values of $x$). Our
computations give the same results as those of \cite{Caulkins 2008}, namely: if
$x_s^1$ and $u_s^1$ and $x_s^2$ and $u_s^2$ are the $x$- and $u$-values at these
steady states, then
\begin{equation}
  \label{eq:steady-states-pequenos}
  \begin{array}{ll}
    x_s^1 \simeq 7.94549\cdot 10^{-7} & x_s^2 \simeq 0.0206096.\\
    u_s^1\simeq 0.0046106 & u_s^2 \simeq 0.284695
  \end{array}
\end{equation}
\begin{figure*}[h!]
  \centering
  \begin{tabular}{ccc}
    \includegraphics{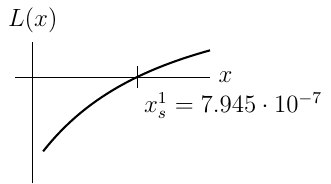}
  % \begin{tikzpicture}[scale=0.6,/tikz/external/up to date check=simple]
  %   \draw[line width=1pt] plot file{tsur-1.table};
  %   \draw (4.5,0) -- (10,0) node[anchor=west]{$x$};
  %   \draw (5,-3) -- (5,1)node[anchor=south]{$L(x)$};
  %   \draw (7.94, 0.3) -- (7.94,-0.3) node at(7.9,-0.8)[anchor=west]
  %   {$x_s^1 = 7.945\cdot 10^{-7}$};
  % \end{tikzpicture}
    &
                     \hspace*{20pt}
                     &
\includegraphics{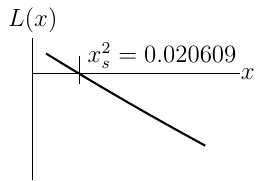}
%                      \begin{tikzpicture}[xscale=1.4,yscale=0.6,/tikz/external/up to date check=simple]
%     \draw[line width=1pt] plot file{tsur-2.table};
%     \draw (1.5,0) -- (4,0) node at(4.1,0){$x$};
%     \draw (1.5,-3) -- (1.5,1.) node at(1.5,1.5){$L(x)$};
%     \draw (2.06,-0.3) -- (2.06,0.5) node[anchor=west] {$x_s^2=0.020609$};
% %    \draw(0,0)--(0.01,0);
% %    \draw (0.0001,0) -- (0.05,0) node[anchor=west]{$x$};
% %    \draw (0.0001,-3) -- (0.0001,3)node[anchor=south]{$L(x)$};
% %    \draw (7.94, 0.3) -- (7.94,-0.3) node at(7.9,-0.5)[anchor=west]{$7.945\cdot 10^{-7}$};
%   \end{tikzpicture}
  \end{tabular}\\
  \text{Figure 6: Stability and instability of steady states $x_s^1$ and $x_s^2$.}
  \label{fig:tsur-1}
\end{figure*}

Tsur's method provides also the nature of these steady states: the first one is
stable (so that solutions tend to it) and the second one is unstable (so that
solutions \emph{never} end on it). These values agree with those of
\cite{Caulkins 2008}. Tsur's method gives the steady states as zeros of the
evolution function $L(x)$. The stability/instability is provided, in Tsur's
method, by the slope of $L(x)$ at those points. Fig. 6 shows how the slope of
$L(x)$ is positive at $x_s^1$, which implies (in a minimization problem) that
$x_s^1$ is stable. On the other hand, the slope of $L(x)$ at $x_s^2$ is
negative, which implies that this is an unstable steady state.

We compute now the dynamic solution approaching the low-value steady state. As
Caulkins \cite{Caulkins 2008} demonstrates, there is a DNS point
$x_D\simeq 0.013$ at which the optimal paths converging towards the high-value
steady state and towards the low-value one have the same cost. We compute our
path starting at $x_D$ in order to verify that our method provides the same
solution as that of \cite{Caulkins 2008}.

Our computations provide, after $6$ seconds, the plots in Figs. 7, 8 and 9 for $x(t)$, $u(t)$ and $\lambda(t)$. After $50$ years we obtain
\begin{equation}
  \label{eq:50-years-low}
  x_s^1 \simeq 7.9445\cdot 10^{-7}, u_s^1 \simeq 0.004612
\end{equation}
and a total cost $\simeq 1.10532$, all in agreement with Caulkins \cite{Caulkins 2008} and the values given by Tsur's algorithm.

\begin{figure*}[h]
  \centering
  \begin{tabular}{ccc}
  \includegraphics{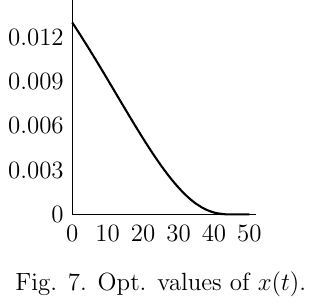}
  % \begin{tikzpicture}[xscale=0.06,yscale=250,/tikz/external/up to date check=simple]
  %   \draw[line width=1pt] plot file{pequenyo-t-x.table};
  %   \draw[anchor=north] node foreach \x in {0,10,20,30,40,50} at (\x,-0.0002) {$\x$};
  %   \draw[anchor=east] node foreach \y in {0.012,0.009,0,0.003,0.006} at (0,\y) {$\y$};
  %   \draw (0,0)--(52,0);
  %   \draw (0,0)--(0,0.0145);
  %   \node at (25,-0.00475) {Fig. 7. Opt. values of $x(t)$.};
  %   % \draw (0,0)--(0,0.32);
  % \end{tikzpicture}
    &
\includegraphics{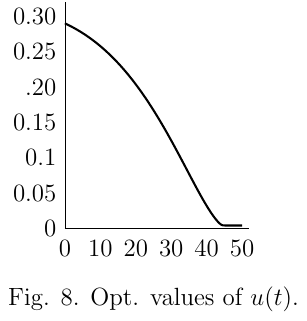}
  %     \begin{tikzpicture}[xscale=0.06,yscale=12,/tikz/external/up to date check=simple]
  %   \draw[line width=1pt] plot file{pequenyo-t-u.table};
  %   \draw[anchor=north] node foreach \x in {0,10,20,30,40,50} at (\x,-0.005) {$\x$};
  %   \draw[anchor=east] node foreach \y in {0,0.05,0.1,0.15,.20,0.25,0.30} at (0,\y) {$\y$};
  %   \draw (0,0)--(52,0);
  %   \draw (0,0)--(0,0.32);
  %   \node at (25,-0.1) {Fig. 8. Opt. values of $u(t)$.};
  % \end{tikzpicture}
    &
      \includegraphics{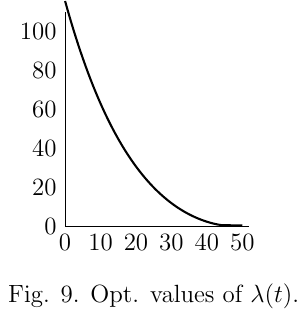}
  % \begin{tikzpicture}[xscale=0.06,yscale=0.033,/tikz/external/up to date check=simple]
  %   \draw[line width=1pt,yshift=36cm] plot file{pequenyo-t-l.table};
  %   \draw[anchor=north,yshift=36cm] node foreach \x in {0,10,20,30,40,50} at (\x,-0.0002) {$\x$};
  %   \draw[anchor=east,yshift=36cm] node foreach \y in {0,20,40,60,80,100} at (0,\y) {$\y$};
  %   \draw[yshift=36cm] (0,0)--(52,0);
  %   \draw[yshift=36cm] (0,0)--(0,110);
  %   \draw[yshift=0cm] (0,0)--(0,0);
  %   \node[] at (25,-0.01) {Fig. 9. Opt. values of $\lambda(t)$.};
  %   % \draw (0,0)--(0,0.32);
  % \end{tikzpicture}
  \end{tabular}
  \label{fig:tsur-2}
\end{figure*}
\pedro{}

\pedro{}\subsection{Parametric study near the High-value steady state}
\pedro{}
After solving the base case, we are going to study the influence of the model
parameters on the steady state solution. In \cite{Caulkins 2008}, the author
presents an extensive sensitivity analysis of each parameter and their
influence both in the state\pedrom{}-control\pedrom{} space and in the state-co-state space. We
pursue a different objective. Using the Algorithm presented in Section 4, we
analyze the influence of the parameters which we deem susceptible of
modification by the government: the efficiency parameters $\gamma$ and $\beta$.

We shall set the steady state of the terrorists' stock as the government
expectation and find the different combinations of $\gamma$ and $\beta$ which
lead to it. First, we perform the sensitivity analysis of each efficiency. The
results are shown in Figs. 10 and 11. We note that in these plots, the cost shown
is not the cost $J$ of the functional over all the time span but the cost by
unit of time of the optimal solution (the integrand). We deem this measure
more adjusted to reality, as once the steady state has been reached, that is
the true value that will set the cost along time.%
\[%
\begin{tabular}
[c]{cc}%
$%
%TCIMACRO{\FRAME{itbpF}{3.1868in}{1.8905in}{0in}{}{}{fig-6.pdf}%
%width 3.1868in;  height 1.8905in;  depth 0in;  original-width 5.1102in;
%original-height 3.0199in;  cropleft "0";  croptop "1";  cropright "1";
%filename 'Fig-6.pdf';file-properties "NPEU";}%
%}}%
%BeginExpansion
{\includegraphics[
height=1.8905in,
width=3.1868in
]%
{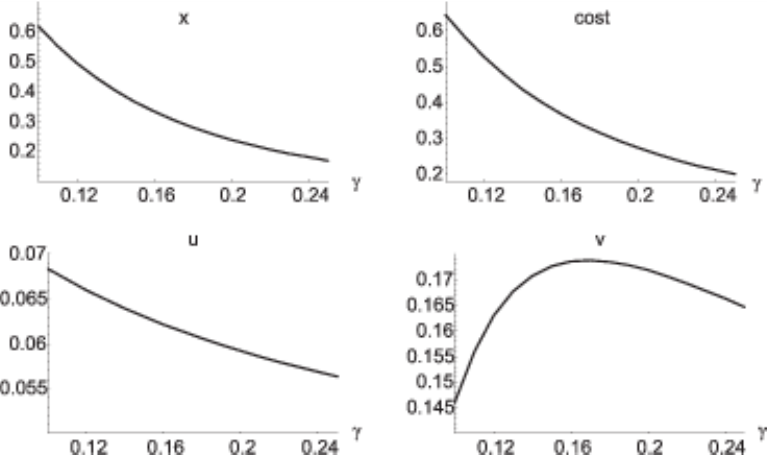}%
}%
%EndExpansion
$ &
%TCIMACRO{\FRAME{itbpF}{3.0528in}{1.8836in}{0in}{}{}{fig-7.pdf}%
%width 3.0528in;  height 1.8836in;  depth 0in;  original-width 5.6818in;
%original-height 3.4938in;  cropleft "0";  croptop "1";  cropright "1";
%filename 'Fig-7.pdf';file-properties "NPEU";}%
%}}%
%BeginExpansion
{\includegraphics[
height=1.8836in,
width=3.0528in
]%
{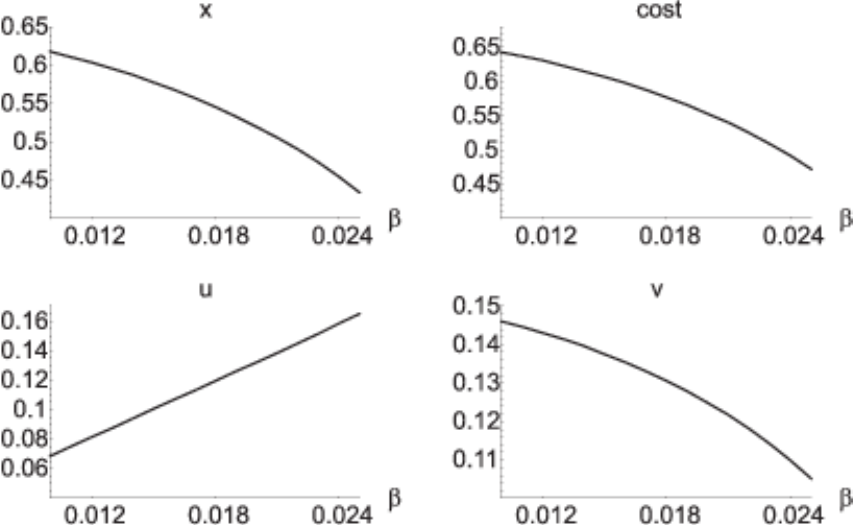}%
}%
%EndExpansion
\\
$\text{Figure 10. Sensitivity analysis of }\gamma\text{.}$ & $\text{Figure 11.
Sensitivity analysis of }\beta\text{.}$%
\end{tabular}
\]
As expected, the influence on the steady state and on cost is greater for
$\gamma$, \pedrom{}as it is related\pedrom{} to the fire actions. The most remarkable fact, in our
view, is that increasing $\beta$ leads to a steady increase of its associated
control $u$ whereas increasing $\gamma$ leads to $v$ reaching a maximum and
then decreasing. The explanation of this lies in the model: for water tactics,
the concave character is included in the model via the $x^{\theta}$ term while
the model for the fire tactics is linear in $x$. Concavity seems to appear
when performing the sensitivity analysis: increasing the efficiency does not
necessarily imply increasing the use of fire tactics.

The joint influence of both efficiencies is studied using successive
applications of our algorithm, using brute force and varying:%
\begin{equation}
\beta\in\lbrack0.01,0.02],\text{ step}=0.001;\text{ \ }\gamma\in
\lbrack0.1,0.2],\text{ step}=0.01
\end{equation}
Our algorithm computes the optimal steady states of the $121$ combinations in
just $33.6\sec$ on a personal computer (Intel Core 2/$2.66$ GHz). This shows
the computational advantage of this method against the one used in the
previous section, which requires the computation of the whole dynamic
solution. Results are shown in Fig. 12.
\[%
\begin{tabular}
[c]{cc}%
$%
%TCIMACRO{\FRAME{itbpF}{2.6368in}{1.785in}{0in}{}{}{fig-8.pdf}%
%width 2.6368in;  height 1.785in;  depth 0in;  original-width 2.9447in;
%original-height 1.983in;  cropleft "0";  croptop "1";  cropright "1";
%filename 'Fig-8.pdf';file-properties "NPEU";}%
%}}%
%BeginExpansion
{\includegraphics[
height=1.785in,
width=2.6368in
]%
{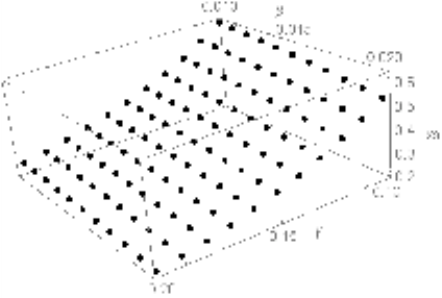}%
}%
%EndExpansion
$ &
%TCIMACRO{\FRAME{itbpF}{1.938in}{1.785in}{0in}{}{}{fig-9.pdf}%
%width 1.938in;  height 1.785in;  depth 0in;  original-width 3.1254in;
%original-height 2.872in;  cropleft "0";  croptop "1";  cropright "1";
%filename 'Fig-9.pdf';file-properties "NPEU";}%
%}}%
%BeginExpansion
{\includegraphics[
height=1.785in,
width=1.938in
]%
{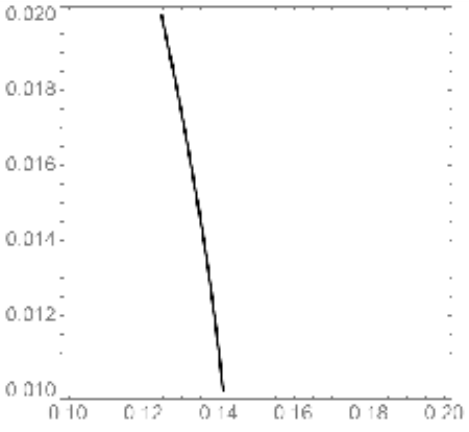}%
}%
%EndExpansion
\\
$\text{Figure 12. Steady-state as function of }\gamma\text{ and }\beta\text{.}$
& $\text{Figure 13. Contour Line }x_{s}(\gamma,\beta)=0.4\text{.}$%
\end{tabular}
\]
Adjusting using least squares, we obtain the function:
\begin{equation}
x_{s}(\gamma,\beta)=1.43186-10.5508\gamma+23.4843\gamma^{2}+1.97376\beta
-301.821\beta^{2}%
\end{equation}
with $r^{2}=0.998$. Using this function we can solve the problem stated at the
beginning: if the government wishes to reach a steady-state of, say,
$x_{s}=0.40$, what are the efficiencies it must get for each tactic? By
imposing%
\begin{equation}
x_{s}(\gamma,\beta)=0.4
\end{equation}
we get the contour curve shown in Fig. 13. Solving the problem for the values
given by this curve gives the results shown in Table 2.
\begin{gather*}
\text{Table 2. Combinations of optimal solutions for }x_{s}(\gamma
,\beta)=0.4\text{.}\\%
\begin{tabular}
[c]{cccccc}%
$\gamma$ & $\beta$ & $x_{s}$ & $u_{s}$ & $v_{s}$ & $cost$\\\hline
0.141 & 0.010 & 0.399 & 0.0638 & 0.1710 & 0.432\\
0.140 & 0.011 & 0.398 & 0.0700 & 0.1693 & 0.432\\
0.139 & 0.012 & 0.398 & 0.0763 & 0.1673 & 0.432\\
0.137 & 0.013 & 0.398 & 0.0826 & 0.1652 & 0.432\\
0.136 & 0.014 & 0.398 & 0.0889 & 0.1628 & 0.432\\
0.134 & 0.015 & 0.397 & 0.0953 & 0.1601 & 0.432\\
0.133 & 0.016 & 0.397 & 0.1017 & 0.1573 & 0.432\\
0.131 & 0.017 & 0.396 & 0.1082 & 0.1541 & 0.431\\
0.129 & 0.018 & 0.395 & 0.1147 & 0.1506 & 0.431\\
0.127 & 0.019 & 0.394 & 0.1213 & 0.1469 & 0.430\\
0.125 & 0.020 & 0.392 & 0.1280 & 0.1427 & 0.429
\end{tabular}
\end{gather*}
Notice first of all that, due to the approximations, $x_{s}$ is not exactly
$0.4$ but the precision is enough. Secondly, observe that modifying the
efficiencies causes (logically) a noticeable change in the associated
controls. However, the most remarkable finding is in the last column: the cost
per time unit of the solution. We see (apart from approximation errors) that
its value is essentially constant on the whole contour curve. This means that
the government has freedom to choose any of the possible combinations so that
the decision should lie on principles other than cost. For example, it may be
easier to increase the efficiency of water strategies rather than fire, due to
ethical or social causes, as these do not require patently bellicose actions.

\section{Sufficient conditions}
We consider now the multi\pedrom{}-control\pedrom{} OCP in its general Bolza form:
\begin{equation}
\underset{u(t)}{\min }J=\int_{0}^{T}F(x(t),u(t),t)dt+B[T,x(T)]  \label{PONT}
\end{equation}%
subject to:
\begin{align}
\dot{x}(t)& =f(x(t),u(t),t);\text{ }x(0)=x_{0} \\
u(t)& \in U(t),\text{ }0\leq t\leq T
\end{align}%
where $x(t)=(x_{1}(t),...,x_{n}(t))\in \mathbb{R}^{n}$ is the state vector and
$u(t)=(u_{1}(t\mathbf),...,u_{n}(t))\in \mathbb{R}^{n}$ the control vector. We assume the following:
\begin{enumerate}[(i)]
\item $F$ and $f = (f_{1}(t),\ldots, f_n(t))$ are continuous.
\item $F$ and $f$ have continuous second derivatives with respect to $t$ and $x$
  but their second derivatives with respect to $u$ may be discontinuous.
\item The control variable $u(t)$ needs only be piecewise continuous.
\item The state variable $x(t)$ is differentiable but its first derivative needs
  only be piecewise continuous (i.e. the graph of $x(t)$ may have ``corners'').
\item The function $B$ has continuous partial derivatives.
\end{enumerate}
The set of admissible controls is, usually, convex and compact.
% Define the Hamiltonian as:
% \begin{equation}
% H(x(t),u(t),\lambda (t),t)=F(x(t),u(t),t)+\lambda (t)f(x(t),u(t),t)
% \end{equation}
% where $\lambda (t)=(\lambda _{1}(t),\ldots ,\lambda _{n}(t))$ is the costate
% vector.
Pontryagin's Theorem establishes necessary conditions for the optimum in
our problem (\ref{PONT}).

In this section, we study the sufficient conditions which guarantee
optimality. They, as is usual
in optimization theory impose some type of convexity of the functions defining
the problem. We concentrate on Arrow's conditions, in this work.

Consider the ``derived Hamiltonian''
\begin{equation}
H(x(t),u(t),\lambda (t),t)=F(x(t),u(t),t)+\lambda (t)f(x(t),u(t),t)
\end{equation}
where $\lambda(t)=(\lambda_1(t),\ldots,\lambda_n(t))$ is the costate vector.
% and define
% \begin{equation*}
% H^{0}(x,\lambda ,t)=\underset{u\in U(t)}{\min }H(x,u,\lambda ,t)
% \end{equation*}
Assume that for the following problem:
\begin{equation}
\underset{u\in U(t)}{\min }H(x,u,\lambda ,t)
\end{equation}
the function:
\begin{equation}
u=u^{0}(x,\lambda ,t)
\end{equation}
is an admissible solution satisfying the PMP. Set:
\begin{equation}
H^{0}(x,\lambda ,t)=H(x,u^{0},\lambda ,t)
\end{equation}
The following result provides a sufficient condition for global optimality, assuming $U$ is compact and convex:

\textbf{Arrow's Theorem}.
\emph{Assume $u^{\ast}(t), x^{\ast}(t)$ and $\lambda^{\ast}(t)$ satisfy PMP for all
$t\in [0,T]$. If $H^0(x, \lambda^{\ast},t)$ is convex in $x$ for each
$t\in [0,T]$ and $B$ is also convex in $x$, then $u^{\ast}$ is an optimal
control for the problem, $x^{\ast}$ is the optimal state trajectory and
$\lambda^{\ast}$ is the optimal trajectory of the costate variables.
}

For local optimality, one only needs to take a compact and convex neighbourhood
of $(x^{\ast}, u^{\ast})$.

\subsection{Infinite-horizon Problems}
If we consider the problem \eqref{PONT} with infinite horizon, following
\cite{Acemoglu 2009}, if the problem includes a \pedrom{}discount factor and\pedrom{} can be stated
as
\begin{equation}
\underset{u(t)}{\min }J=\int_{0}^{\infty }F(x(t),u(t))e^{\pedrom{}-r\pedrom{} t}dt
\end{equation}
subject to:
\begin{align}
\dot{x}(t)& =f(x(t),u(t),t);\text{ }x(0)=x_{0} \\
u(t)& \in U(t)
\end{align}
then:

\textbf{Sufficient conditions for an infinite-horizon problem with discount
  factor}.  \emph{Given the problem above with the same conditions on $F$, $f$,
  $x$ and $u$ and with $H$ and $H^0$ defined as before, let
  $u^{\ast}(t), x^{\ast}(t)$ and $\lambda^{\ast}(t)$ be admissible solutions
  satisfying PMP for all $t$. Then, if
  \begin{equation}
    \lim_{t\rightarrow \infty }\lambda ^{\ast }(t)x(t)\geq 0
  \end{equation}
  for all admissible $x(t)$ and the conditions for Arrow's Theorem hold, then
  the pair $(x^{\ast}(t), u^{\ast}(t))$ is a global minimum.  }

Of course, the same principle as above serves to ensure a local minimum. In our
case, the first inequality
  \begin{equation}
    \lim_{t\rightarrow \infty }\lambda ^{\ast }(t)x(t)\geq 0
  \end{equation}
holds by the conditions imposed on the problem.

% N\'{o}tese que en el teorema anterior se han fusionado las condiciones de
% concavidad de Mangasarian en las funciones $F$ y $f$ y la condici\'{o}n de
% no negatividad en $\lambda (t)$ en una \'{u}nica condici\'{o}n de concavidad
% en el hamiltoniano $H$. La concavidad de $H$ es en $(x,u)$, de manera
% conjunta. Por el contrario, la condici\'{o}n de Arrow es que $H^{0}$ sea c
% \'{o}ncava en la variable $x$ solamente. Las dos primeras condiciones ya las
% hemos analizado antes.

% En cuanto a la \'{u}ltima:
% \begin{equation}
% \lim_{t\rightarrow \infty }\lambda ^{\ast }(t)x(t)\geq 0
% \end{equation}
% su cumplimiento tambi\'{e}n es evidente, por las hip\'{o}tesis del problema.

% Por tanto las conclusiones son id\'{e}nticas para ambos problemas.\newpage

\subsection{A new functional}
Some computations show that we cannot guarantee Arrow's conditions with the
previous functional for a solution tending to $x_s$. Despite this, the value of
the corresponding second partial derivative is near enough to zero to prevent us
from also stating that it is \emph{not} an optimum.

Recall that the functional is given by:
\begin{equation}
  (cx+u^{2}+v^{2})e^{-rt}
\end{equation}
so that the cost is considered linear in the number of terrorists. This is an
assumption which we think does not take into account the value of information
sharing, interaction and synergies in groups of people. In our view, a convex
function of $x$ provides an improved cost model in terms of the number of
terrorists. We are going to perform the study solving (for simplicity we take
$x^2$) the new problem:
\begin{align}
  \underset{u,v\geq 0}{\min }J &=
   \underset{u,v\geq 0}{\min }\int_{0}^{\infty}
    (cx^{2}+u^{2}+v^{2})e^{-rt}dt  \label{OCP}\\
  \dot{x} &=
  \tau +(1+\rho v)kx^{\alpha }-
    \mu x-\beta \ln (1+u)x^{\theta}-
    \gamma \ln (1+v)x;\text{ \ }x(0)=x_{0}  \\
  u(t) &\geq 0;\text{ \ }v(t)\geq 0
\end{align}

With the same data as in the previous case, the solution we find leads to a
high-value steady-state, somewhat different from the one obtained before. In
this case, we have:
\begin{align}
x_{s}& =0.605 \\
u_{s}& =0.081 \\
v_{s}& =0.163
\end{align}
The Hamiltonian is given by:
\begin{equation}
  \begin{split}
    H(x,u,v,\lambda ,t)&=
    e^{-rt}\left( cx^{2}+u^{2}+v^{2}\right) \\
    &+\lambda \left( \tau +(1+\rho
      v)kx^{\alpha }-\mu x-\beta \ln (1+u)x^{\theta }-\gamma \ln (1+v)x\right)
\end{split}
\end{equation}
so that:
\begin{align}
H_{u}(x,u,v,\lambda ,t) &=2ue^{-rt}-\frac{\beta \lambda x^{\theta }}{u+1} \\
H_{v}(x,u,v,\lambda ,t) &=2ve^{-rt}+\lambda \left( k\rho x^{\alpha }-\frac{
\gamma x}{v+1}\right)
\end{align}
The Hamiltonian is convex in $u$ and $v$, as:
\begin{align}
H_{{u}u}(x,u,v,\lambda ,t) &=2e^{-rt}+\frac{\beta \lambda x^{\theta }}{
(u+1)^{2}}\geq 0 \\
H_{vv}(x,u,v,\lambda ,t) &=2e^{\pedrom{}-rt\pedrom{}}+\frac{\gamma \lambda x}{
(v+1)^{2}}\geq 0
\end{align}
Solving the equations:
\begin{align}
H_{u}(x,u,v,\lambda ,t) &=0 \\
H_{v}(x,u,v,\lambda ,t) &=0
\end{align}
we obtain:
\begin{align}
&u^{0}(x,\lambda ,t) =\frac{1}{2}\left( \sqrt{2\beta \lambda
    e^{rt}x^{\theta }+1}-1\right)\\
&v^{0}(x,\lambda ,t) =\frac{1}{4}\left( \sqrt{k^{2}\lambda ^{2}\rho
^{2}e^{2rt}x^{2\alpha }-4k\lambda \rho e^{rt}x^{\alpha }+8\gamma \lambda
xe^{rt}+4}+k\lambda \rho \left( -e^{rt}\right) x^{\alpha }-2\right)
\end{align}
so that the derived Hamiltonian is given by:
\begin{equation}
  \begin{split}
    H^{0}(x,\lambda ,t)= &\; e^{-rt}\left( cx^{2}+u^0(x,\lambda,t)^2+
      v^0(x, \lambda, t)^2 \right)\\
    &+\lambda \left( \tau +\left( 1+\rho v^0(x, \lambda, t) \right) kx^{\alpha
      }-\mu x-\beta
      \ln \left(1 + u^0(x, \lambda, t) \right) x^{\theta} \right.\\
    & - \lambda \gamma \ln \left( 1+v^0(x, \lambda, t) \right) x
\end{split}
\end{equation}
For the $\lambda^{\ast}(t)$ obtained, we get the plot of $H^0_{xx}$ as a
function of $t$ shown in Fig. 14, as a function of $t$, which guarantees the
local convexity and hence, local optimality of our solution:

\begin{gather*}%
%TCIMACRO{\FRAME{itbpF}{2.2866in}{1.075in}{0in}{}{}{fig-2.pdf}%
%width 2.2866in;  height 1.075in;  depth 0in;  original-width 6.0727in;
%original-height 2.8392in;  cropleft "0";  croptop "1";  cropright "1";
%filename 'Fig-2.pdf';file-properties "NPEU";}%
%}}%
%BeginExpansion
  {\includegraphics
{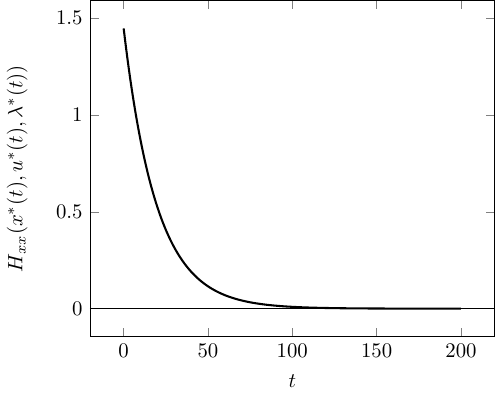}%
}%
%EndExpansion
\\
\text{Figure 14. Convexity of the derived Hamiltonian.}%
\end{gather*}

\section{Conclusions}

We consider ways for a government to optimally employ \textquotedblleft
water\textquotedblright\ and \textquotedblleft fire\textquotedblright%
\ strategies for fighting terrorism. The model tries to balance the costs of
terror attacks with the cost of terror control. We present two main
contributions. First, a new effective algorithm for computing the dynamical
solution whose cyclic nature allows its use in models of greater dimension
(say with more controls or more state variables) without any conceptual
modification. Secondly, we present a method for computing the steady-state
solution which can be used to analyze the dependence of the steady-state
strategy on several parameters.

We show how this latter method allows us to compute, a priori, the optimal
values of the steady-state, without having to solve the whole dynamical
problem. This way one can analyze the influence of the parameters in a simpler
way, which leads to remarkable computational savings. We can study not just
the influence of the parameters on the model but (what is more important),
compute what the parameters must be in order to obtain a desired steady-state.
We have studied the influence of the efficiency parameters $\gamma$ and
$\beta$. We consider that the study of the influence of the other parameters
may be interesting from the theoretical point of view but, as they are
intrinsic to the model, they cannot be modified (at least in a simple way) by
the authorities. However, the efficiency $\gamma$ of the fire controls is a
vital question for any army, as is the efficiency $\beta$ of the water
controls for any government. In the latter case, both the secret services
intelligence and even the educational measures taken by the governments are relevant.

Finally, we study a modification of the model in which the cost function is
strictly convex in the number of terrorists and verify that the solution found
by our method is, indeed, a local optimum.


\begin{thebibliography}{99}                                                                                               %
\bibitem {Heymann 2003}P. Heymann, Dealing with terrorism after september 11,
2001: An overview. In Countering terrorism: Dimensions of preparedness, MIT
Press, 57-72, 2003.

\bibitem {Keohane 2003}N.O. Keohane, R.J. Zeckhauser, The Ecology of Terror
Defense. Journal of Risk and Uncertainty 26, 201-229, 2003.

\bibitem {Castillo 2003}C. Castillo-Chavez, B. Song, Models for the
transmission dynamics of fanatical behaviors. In Bioterrorism: Mathematical
Modeling Applications in Homeland Security, SIAM, Philadelphia, 155-172, 2003.

\bibitem {Kaplan 2005}E.H. Kaplan, A. Mintz, S. Mishal, C. Samban, What
happened to suicide bombings in Israel? Insights from a terror stock model.
Studies in Conflict and Terrorism 28, 225-235, 2005.

%\bibitem {Abadie 2003}A. Abadie, J. Gardezabal, The economic costs of
%conflict: a case control study for the Basque country. American Economic
%Review 93, 113-132, 2003.

%\bibitem {Barros 2006}C.P. Barros, J. Passos, L. Gil-Alana, The timing of the
%ETA terrorist attacks. Journal of Policy Modeling 28, 335-346, 2006.

\bibitem {Udwadia 2006}F. Udwadia, G. Leitmann, L. Lambertini, A dynamical
model of terrorism. Discrete Dynamics in Nature and Society, 1-32, 2006.

\bibitem {Caulkins 2008}J.P. Caulkins, D. Grass, G. Feichtinger, G. Tagler,
Optimizing counter-terror operations: Should one fight fire with
\textquotedblleft fire\textquotedblright\ or \textquotedblleft
water\textquotedblright?. Computers \& Operations Research 35, 1874-1885, 2008.

\bibitem {Kress 2009}M. Kress, R. Szechtman, Why Defeating Insurgencies Is
Hard: The Effect of Intelligence in Counterinsurgency Operations-A Best-Case
Scenario, Operations Research 57(3), 578-585, 2009.

\bibitem {Faria 2012}J.R. Faria, D. Arce, Counterterrorism And Its Impact On
Terror Support And Recruitment: Accounting For Backlash. Defence and Peace
Economics 23(5), 431-445, 2012.

\bibitem {Seidel 2016}A. Seidl, E.H. Kaplan, J.P. Caulkins, S. Wrzaczek, G.
Feichtinger, Optimal control of a terror queue, European Journal of
Operational Research 248, 246-256, 2016.

\bibitem {Grass 2008}D. Grass, J.P. Caulkins, G. Feichtinger, G. Tragler, D.A.
Behrens, Optimal Control of Nonlinear Processes: With Applications in Drugs,
Corruption and Terror. Springer-Verlag, Berlin, 2008.

\bibitem {Schneider 2015}F. Schneider, T. Br\"{u}ck, D. Meierrieks, The
Economics of Counter-Terrorism: A Survey, Journal of Economic Surveys 29(1),
131-157, 2015.

\bibitem {Tsur 2001}Y. Tsur and A. Zemel, The infinite horizon dynamic
optimization problem revisited: A simple method to determine equilibrium
states, European Journal of Operational Research 131(3), 482-490, 2001.

\bibitem {Luo 1992}Z.Q. Luo and P. Tseng, On the convergence of the coordinate
descent method for convex differentiable minimization, J. Optim. Theory Appl.
72(1), 7-35, 1992.

\bibitem {Zangwil 1969}W.L. Zangwill, Nonlinear Programming: A Unified
Approach, Prentice Hall, Nueva Jersey, 1969.

\bibitem {Pontryagin 1987}L.S. Pontryagin, Mathematical Theory of Optimal
Processes (Classics of Soviet Mathematics), CRC Press, 1987.

\bibitem {Clarke 1983}F.H. Clarke, Optimization and Nonsmooth Analysis,
Wiley-Interscience, New York, 1983.

\bibitem {Chiang 2000}A. Chiang, Elements of Dynamic Optimization. Waveland
Press, 2000.

\bibitem {Benveniste 1982}L.M. Benveniste, J.A. Scheinkman, Duality Theory for
Dynamic Optimization Models of Economics: The Continuous Time Case, Journal of
Economic Theory 27, 1-19, 1982.

\pedrom{}\bibitem{aseev}S.M. Aseev, A.V. Kryazhimskiy, The Pontryagin Maximum Principle and Transversality Conditions for a Class of Optimal Control Problems with Infinite Time Horizons, SIAM Journal on Control and Optimization 43(3), 1094-1119, 2004.

\bibitem{bolza}L. Bay\'on, J.M. Grau, M.M. Ruiz, P.M. Su\'arez, A Bolza Problem in Hydrothermal Optimization, Applied Mathematics and Computation 184(1), 12-22, 2007.
\pedrom{}

\bibitem {Corless 1996}R. M. Corless, G. H. Gonnet, D. E. G. Hare, D. J.
Jeffrey and D. E. Knuth, On the Lambert W Function, Adv. Comput. Math.
\textbf{5}, 329-359 (1996)

\bibitem {Acemoglu 2009}D. Acemoglu, Introduction to Modern Economic Growth,
  Princeton University Press, 2009.
\end{thebibliography}
\end{document}